\documentclass{amsart}

\usepackage{amsmath}
\usepackage{amscd}
\usepackage{amssymb} 

\newcommand{\cal}{\mathcal}
\newcommand{\bk}{{\bf k}}

\DeclareMathOperator{\Img}{Im}

\DeclareMathOperator{\Ker}{Ker}

\newtheorem{theorem}{Theorem}[section]
\newtheorem{theorem/definition}{Theorem/Definition}[section]
\newtheorem{proposition}{Proposition}[section]
\newtheorem{lemma}{Lemma}[section]

\theoremstyle{remark}
\newtheorem{remark}{Remark}[section]

\theoremstyle{definition}

\begin{document}

\title
{Identification of Two Frobenius Manifolds In Mirror Symmetry}
\author{Huai-Dong Cao \& Jian Zhou}
\address{Department of Mathematics\\
Texas A \& M University\\
College Station, TX 77843}
\email{cao@math.tamu.edu,  zhou@math.tamu.edu}
\begin{abstract}
We identify two Frobenius manifolds obtained from
two different differential 
Gerstenhaber-Batalin-Vilkovisky algebras on a compact K\"ahler manifold.
One is constructed on the Dolbeault  cohomology in Cao-Zhou \cite{Cao-Zho2}, 
and the other on the de Rham 
cohomology in the present paper. 
This can be considered as a generalization of the identification of 
the Dolbeault cohomology ring 
with the complexified de Rham cohomology ring on a K\"ahler manifold. 
\end{abstract}
\maketitle
\date{}

\footnotetext[1]{Both authors were supported in part by NSF }

The mysterious Mirror Conjecture \cite{Yau} in string theory has 
enabled the physicists to write down a formula 
\cite{Can-Oss-Gre-Park} 
on the number of rational curves of any degree 
on a quintic in ${\Bbb C}{\Bbb P}_4$. 
Recently, Lian-Liu-Yau \cite{Lia-Liu-Yau} have given 
a rigorous proof of this important formula.

The  rapidly progressing theory of quantum cohomology, 
also suggested by physicists,
has lead to a better mathematical formulation of
the Mirror Conjecture.
Now a version of the Mirror Conjecture can be formulated  as 
the identification of Frobenius manifold structures obtained by
different constructions.
(For an exposition of this point of view, 
the reader is referred to a recent paper by Manin \cite{Man}.)
More precisely, on a Calabi-Yau manifold $X$, 
there are two natural algebras 
\begin{eqnarray*}
A(X) = \oplus_{p, q} H^q(X, \Omega^p), &
B(X) = \oplus_{p, q} H^q(X, \Omega^{-p}),
\end{eqnarray*}
where  $\Omega^{-p}$ is the sheaf of holomorphic sections
to $\Lambda^pTX$.
By Bogomolov-Tian-Todorov theorem, 
the moduli space of complex structures on $X$ 
is an open subset in $H^1(X, \Omega^{-1})$.
Witten \cite{Wit} suggested the construction of 
an extended moduli space of complex structures.
In the case of Calabi-Yau manifold,
Barannikov-Kontsevich \cite{Bar-Kon} generalized the construction 
of Tian \cite{Tia} and Todorov \cite{Tod}
to show that the extended moduli space is a supermanifold with Bosonic part 
an open set in 
$$B^{\text{even}}(X) = \oplus_{p+q=\text{even}} 
	H^q(X, \Omega^{-p}).$$
Furthermore, the same method is essential to construction of 
a structure of Frobenius supermanifold 
on this extended moduli space given by 
Barannikov-Kontsevich \cite{Bar-Kon}.
This construction of Frobenius manifold structure 
has been generalized by
Manin \cite{Man} to general differential 
Gerstenhaber-Batalin-Vilkovisky (dGBV)
algebras with some mild conditions.  
On the other hand, by  Hodge theory,
$A(X)$ is isomorphic to the de Rham cohomology 
with complex coefficients $H^*_{dR}(X, {\Bbb C})$.
Now  the complexified deformations  of K\"{a}hler 
form is parameterized by 
$H^{1, 1}(X) \subset H^*_{dR}(X, {\Bbb C})$.
On the other hand, there is a construction of Frobenius manifold structure on
the de Rham cohomology, 
which was rigorously established in the work of Ruan-Tian \cite{Rua-Tia}
on the mathematical formulation of quantaum cohomology. 
For a survey of quantum cohomology, see Tian \cite{Tia}. 
For some remarks on some recent development, see the introduction of
Li-Tian \cite{Li-Tia}.
In Cao-Zhou \cite{Cao-Zho2}, the authors  gave a construction of
Frobenius supermanifold modelled on $A(X)$, and conjectured
it can be identified with the construction via quantum cohomology.
Notice that in general there is a problem of convergence in
the construction via quantum cohomology.
This problem has been solved only 
in the case of complete intersection Calabi-Yau manifolds.
See Tian \cite{Tia}.
In our construction, a standard argument in 
Kodaira-Spencer-Kuranishi 
deformation theory \cite{Mor-Kod} guarantees the convergence .

In this paper, we give  a construction of 
a Frobenius supermanifold modelled on $H^*_{dR}(X, {\Bbb C})$.
We then identify it with the Frobenius supermanifold modelled 
on $A(X)$ constructed in Cao-Zhou \cite{Cao-Zho2}.
As to the knowledge of the authors, 
this gives the  first nontrivial result on the identification of
two Frobenius  supermanifolds.

{\bf Acknowledgements}. 
{\em The authors would like to thank Gang Tian and S.-T. Yau for 
their interest. The work in this paper is carried out
while the second author is visiting Texas A$\&$M University.
He likes to express his appreciation for the hospitality
and financial support of the Mathematics Department 
and the Geometry-Analysis-Topology group.}

\section{A construction of Frobenius (super)manifolds}

In this section, we review a construction of Frobenius supermanifolds.
For details, the reader should consult the papers
by Tian \cite{Tia}, Todorov \cite{Tod}, 
Barannikov-Kontsevich \cite{Bar-Kon} and Manin \cite{Man}.
Here, we follow the formulation in Manin \cite{Man}.

\subsection{Gerstenhaber algebras}
A {\em Gerstenhaber algebra} consists of a triple 
$({\cal A} = \oplus_{i \in {\Bbb Z}_2} {\cal A}^i, 
\wedge, [\cdot \bullet \cdot])$,
such that $({\cal A}, \wedge)$ is a ${\Bbb Z}_2$-graded 
commutative associative algebra over a filed ${\bk}$, 
$({\cal A}[-1], [\cdot \bullet \cdot])$ is a graded Poisson algebra
with respect to the multiplication $\wedge$.
Here ${\cal A}[-1]$ stands for the vector space ${\cal A}$ with a
new grading: ${\cal A}[-1]^i ={\cal A}^{i+1}$.

An operator $\Delta$ of odd degree is said to generate 
the Gerstenhaber bracket if for all homogeneous $a, b \in {\cal A}$,
$$[a \bullet b] = (-1)^{|a|}(\Delta(a \wedge b) - \Delta a \wedge b 
- (-1)^{|a|} a \wedge \Delta b).$$
To the authors' knowledge,
a formula of this type in differential geometry was first discovered 
by Tian \cite{Tia1} (see also Todorov \cite{Tod}) to prove the important
result that deformations of Calabi-Yau manifolds are unobstructed.

\subsection{Gerstenhaber-Batalin-Vilkovisky algebras}
If there is an operator $\Delta$ on a Gerstenhaber algebra ${\cal A}$
which generates the Gerstenhaber bracket,
such that $\Delta^2 = 0$, 
then ${\cal A}$ is called a 
{\em Gerstenhaber-Batalin-Vilkovisky algebra} (GBV algebra).

\subsection{Differential Gerstenhaber-Batalin-Vilkovisky algebras}
A  {\em differential Gerstenhaber-Batalin-Vilkovisky algebra}
(dGBV algebra) is 
a GBV algebra with a $\bk$-linear derivation $\delta$ of odd degree
with respect to $\wedge$, such that
$$ \delta^2 = \delta\Delta + \Delta \delta = 0.$$
We will be interested in the cohomology group $H({\cal A}, \delta)$.

\subsection{Integral on dGBV algebras}
A $\bk$-linear functional $\int: {\cal A} \rightarrow \bk$ 
on a dGBV-algebra is called {\em an integral} if
for all $a, b \in {\cal A}$,
\begin{eqnarray}
\int (\delta a) \wedge b & = & 
	(-1)^{|a|+1} \int a \wedge \delta b, \label{int1} \\
\int (\Delta a ) \wedge b & = &
	(-1)^{|a|} \int a \wedge \Delta b. \label{int2}
\end{eqnarray} 
Under these conditions, it is clear that $\int $ induces a
scalar product on $H = H({\cal A}, \delta)$:
$(a, b ) = \int a \wedge b$.
If it is nondegenerate on $H$, we say that the integral
is {\em nice}.
It is obvious that 
$$(\alpha \wedge \beta, \gamma) = (\alpha, \beta \wedge \gamma).$$
By definition, $(H({\cal A}, d), \wedge, (\cdot, \cdot))$ is then
a {\em Frobenius algebra}, when ${\cal A}$ has a nice integral.

\subsection{Formal Frobenius supermanifolds from dGBV algebras}

To obtain a deformation of the ring structure on $H$,
consider $\delta_a: {\cal A} \rightarrow {\cal A}$ for even 
$a \in {\cal A}$.
When $a$ satisfies
\begin{eqnarray*}
&\delta a + \frac{1}{2}[a \bullet a] = 0, \\
&\Delta a = 0, 
\end{eqnarray*}
then $({\cal A}, \wedge, \delta_a, 
\Delta, [\cdot \bullet \cdot])$ is also a dGBV algebra.
If $\int$ is an integral for 
$({\cal A}, \wedge, \delta, \Delta, [\cdot \bullet \cdot])$, 
so is it for 
$({\cal A}, \wedge, \delta_a, \Delta, [\cdot \bullet \cdot])$. 

Let $K$ be the formal super power series generated by the
super vector space 
$H = H^{\text{even}} \oplus  H^{\text{odd}}$. 
Let $\{e_j, j =0, \cdots \}$ be 
a basis of $H$ consist of homogeneous elements,
such that $e_0 = 1$, 
$\{ x^j: j = 0, \cdots  \}$ be the dual basis of $H'$.
Then $K$ is generated by $x^j$, with the relations
$$ x^j x^k = (-1)^{|x^j||x^k|}x^kx^j. $$
Hence $K$ is isomorphic to $S^*(H^{\text{even}}) \otimes 
\Lambda^*(H^{\text{odd}})$.
Let ${\cal A}_K = K \otimes {\cal A}$,
and extend $\wedge$, $\delta$, $\Delta$ 
and $[\cdot\bullet\cdot]$ as super-(bi)module maps.
The construction of Frobenius manifold structure is based on the existence of a solution 
$\Gamma = \sum \Gamma_n$ to
\begin{eqnarray*}
&\delta \Gamma + \frac{1}{2}[\Gamma \bullet \Gamma] = 0, \\
&\Delta \Gamma = 0, 
\end{eqnarray*}
such that $\Gamma_0 = 0$,
$\Gamma_1 = \sum x^j e_j$, $e_j \in \Ker \delta \cap \Ker \Delta$.
For $n > 1$, $\Gamma_n \in \Img \Delta$ is 
a homogeneous super polynomial of degree $n$ in $x^j$'s,
such that the total degree of $\Gamma_n$ is even. 
Furthermore, $x^0$ only appears in $\Gamma_1$.
Such a solution is called a {\em normalized universal solution}.
Its existence can be established inductively. 
This is how Tian \cite{Tia1} and Todorov \cite{Tod} proved
that the deformation of complex structures on a Calabi-Yau
manifold is unobstructed.
It was later generalized by Barannikov-Kontsevich \cite{Bar-Kon}
to the case of extended moduli space of  complex structures
of a Calabi-Yau manifold.
Manin \cite{Man} further generalized it to the case of 
dGBV algebras.
In the rest of this section,  
we briefly review some relevant materials 
$\S 4 - \S 6$ in Manin \cite{Man}.

\begin{lemma} \label{keylemma}
Let $({\cal A}, \wedge, \delta, \Delta, [\cdot\bullet\cdot])$  be
a dGBV algebra.

(a) The following statements are  equivalent:

(i) The inclusions $i: (\Ker \Delta, \delta) \hookrightarrow
({\cal A}, \delta)$
and $j: (\Ker \delta, \Delta) \hookrightarrow
({\cal A}, \Delta)$
induce isomorphisms on cohomology.

(ii) We have the following equalities:
\begin{align*}
\Img \delta\Delta = \Img \Delta \delta = \Img \delta \cap \Ker \Delta, \tag{A} \\
\Img \delta\Delta = \Img \Delta \delta = \Img \Delta \cap \Ker \delta \tag{B}.
\end{align*}

(iii) We have the following equality:
\begin{align*}
\Img \delta\Delta = \Img \Delta \delta = 
(\Ker \delta \cap \Ker \Delta) \cap (\Img \delta + \Img \Delta).
\tag{C}
\end{align*}

(b) If any one of the statements in (a) holds, 
then we have isomorphisms
$$H({\cal A}, \delta) \cong H({\cal A}, \Delta) \cong 
(Ker \delta \cap \Ker \Delta) / \Img \delta\Delta.$$

(c) If any one of the statements in (a) holds, 
then the quotient map 
$\Ker \Delta \rightarrow H({\cal A}, \Delta)$
is a surjective homomorphism of 
differential graded odd Lie algebras
$$\phi: (\Ker \Delta, \delta, [\cdot\bullet\cdot])
\rightarrow (H({\cal A}, \Delta), 0, 0),$$
which induces isomorphism on cohomology.

(d) If any one of the statements in (a) holds, 
then there exists a normalized universal solution.
\end{lemma}

\begin{proof}
(a.i) and (a.ii) and (b) are part of the more general
statements of Lemma 5.4.1 in Manin \cite{Man},
(a.iii) is part of Remark 5.4.2  in \cite{Man}. 
The proofs of (c) and (d) are given in the proofs
of Proposition 6.1.1 in \cite{Man}, 
which in turn uses  Theorem 4.2 in \cite{Man}.
\end{proof}

Given a normalized universal solution $\Gamma$,
we have a universal shifted dGBV  algebra
$({\cal A}_{\Gamma}, \wedge, \delta_{\Gamma}, \Delta,
[\cdot\bullet\cdot])$.
There is a canonical way to  extend a $\delta$-closed
element to a $\delta_{\Gamma}$-closed element.
For any $X \in \Ker \delta \cap \Ker \Delta$,
let $\overline{X}: K \rightarrow K$ be the derivation
given by supercontraction by the class $[X]$ of $X$ 
in $H({\cal A}, \delta)$.
Tensor with $\Ker \Delta$, we obtain a right $\Ker \Delta$-linear
derivation, ${\cal A}_K \rightarrow {\cal A}_K$,
still denoted by $\overline{X}$,
since there is no danger of confusion.
Now the leading term of $\overline{X}\Gamma$ is $X$,
and 
$\overline{X}(\delta \Gamma + \frac{1}{2}[\Gamma, \Gamma]) = 0$
is equivalent to $\delta_{\Gamma} (\overline{X}\Gamma) = 0$.
So we get a canonical way to extend $X$ to a 
$\delta_{\Gamma}$-closed element.
Conversely, any element of 
$H_{\Gamma} =H({\cal A}_\Gamma, \delta_{\Gamma})$
with leading term $X$ is obtained this way.
See Corollary 4.2.1 in Manin \cite{Man}.
Apply this result, one  can define a multiplication $\circ$
on $H_{\Gamma}$ as follows:
$$\overline{X \circ Y} \Gamma \equiv (\overline{X}\Gamma) \wedge
(\overline{Y}\Gamma) \pmod{\Img\delta_{\Gamma}}.$$

Assume now that  $\int$ is an integral 
for $({\cal A}, \wedge, \delta, \Delta,[\cdot \bullet \cdot])$, 
so is it for 
$({\cal A}, \wedge, 
\delta_{\Gamma}, \Delta, [\cdot \bullet \cdot])$. 
For $X, Y \in \Ker \delta \cap \Ker  \Delta$,
we have 
\begin{eqnarray*}
\overline{X}\Gamma = X + \sum_{m \geq 2} \overline{X}\Gamma_m, &
\overline{Y}\Gamma = Y + \sum_{n \geq 2} \overline{Y}\Gamma_n.
\end{eqnarray*}
Now $\overline{X}\Gamma_m$ and $\overline{Y}\Gamma_n$
lies in $\Img \Delta$, for $m \geq 2$ and $n \geq 2$,
by the fact that 
\begin{eqnarray*}
\int a \wedge b = 0, & 
\text{if} \; \; a \in \Ker \Delta, \; b \in 
\Img \Delta ,
\end{eqnarray*}
we have
\begin{eqnarray*}
&   & \int \overline{X}\Gamma \wedge \overline{Y}\Gamma \\
& = & \int X \wedge Y 
+ \int X \wedge \sum_{n \geq 2} \overline{Y}\Gamma_n
+ \int \sum_{m \geq 2} \overline{X}\Gamma_m \wedge Y 
+ \int \sum_{m \geq 2} \overline{X}\Gamma_m \wedge 
\sum_{n \geq 2} \overline{Y}\Gamma_n \\
& = & \int X \wedge Y.
\end{eqnarray*}
(Cf. Barannikov-Kontsevich \cite{Bar-Kon}, Claim 6.2.)
Now we can state the main result of Manin \cite{Man}, $\S 6$.

\begin{theorem} \label{thm:constrution}
Let $({\cal A}, \wedge, \delta, \Delta, [\cdot \bullet \cdot])$
be a dGBV algebra satisfying the following conditions:
\begin{enumerate}
\item $H = H({\cal A}, \delta)$ is finite dimensional.
\item There is a nice integral on ${\cal A}$.
\item The inclusions $(Ker \Delta, \delta) \hookrightarrow 
({\cal A}, \delta)$ and $(Ker \delta, \Delta) \hookrightarrow 
({\cal A}, \Delta)$  induce isomorphisms of cohomology.
\end{enumerate}
Then there is a structure of a formal Frobenius manifold on 
the formal spectrum of $K$, 
the algebra of formal superpower series generated by $H'$,
where $H'$ is the dual $\bk$-vector space of $H$.
\end{theorem}

\section{Frobenius supermanifolds modelled on de Rham cohomology}

\subsection{dGBV algebra in Poisson geometry}
Our reference  for Poisson geometry is Vaisman \cite{Vai}.
For Gerstenhaber algebra and dGBV algebra in Poisson geometry,
see Xu \cite{Xu} and the references therein.

Let $w \in \Gamma(X, \Lambda^2TX)$ be a bi-vector field.
It is called a Poisson bi-vector if the Schouten-Nijenhuis
bracket of $[w, w]$ vanishes \cite{Vai}.
Let ${\cal A} = \Omega(X)$ with the ordinary wedge product $\wedge$,
and the exterior differential $d$.
Following Koszul \cite{Kos}, 
define $\Delta: \Omega^*(X) \rightarrow \Omega^{*-1}(X)$
by $\Delta \alpha = w \vdash d \alpha - d (w \vdash \alpha)$,
for $\alpha \in \Omega^*(X)$, where $\vdash$ is the contraction.
Koszul \cite{Kos} proved that 
$\Delta^2 = 0$ and $d \Delta + \Delta d = 0$.
Also defined in Koszul \cite{Kos} is the {\em covariant 
Schouten-Nijenhuis bracket}
$$\{\lambda, \mu \} = (\Delta \lambda) \wedge \mu 
+ (-1)^{|\lambda|} \lambda \wedge (\Delta \mu) 
- \Delta ( \lambda \wedge \mu).$$
This bracket has the following properties
(see e.g. Proposition 4.24 in Vaisman \cite{Vai}, 
where $\delta$ is used instead of $\Delta$):
\begin{eqnarray*}
&\{ \lambda, \mu \} = (-1)^{|\lambda||\mu|}
	\{ \mu, \lambda \}, \\
&(-1)^{|\lambda|(|\nu|-1)}\{\lambda, \{\mu, \nu\}\} 
+(-1)^{|\mu|(|\lambda|-1)} \{\mu, \{\nu, \lambda\}\}
+(-1)^{|\nu|(|\mu|-1)} \{\nu, \{\lambda, \mu\}\} = 0,\\
&\{\lambda, \mu \wedge \nu \} = \{ \lambda, \mu\} \wedge \nu
+ (-1)^{|\lambda|(|\mu|+1)} \mu \wedge \{\lambda, \nu\}.
\end{eqnarray*}
Now it is clear that if we define 
$$[\alpha \bullet \beta]  = (-1)^{|\alpha|-1}\{\alpha, \beta\},$$
then $({\cal A}, \wedge, \delta = d, \Delta, [\cdot \bullet \cdot])$ is
a dGBV algebra. 
Such a dGBV algebra  is known in Poisson geometry (see Xu \cite{Xu}).
The authors are led to this by our work on
quantum de Rham cohomology \cite{Cao-Zho1}.

When $X$ is closed, let $\int: {\cal A} \rightarrow {\Bbb R}$
be the ordinary integral of differential forms over $X$.
Then clearly $(1)$ is satisfied.
To check $(2)$, we need the following

\begin{lemma} \label{lem:contraction}
If $\alpha, \beta \in \Omega^*(X)$ satisfy 
$|\alpha| + |\beta| = \dim (X) + 2$,
then we have
$$\int_X (w \vdash \alpha) \wedge \beta = 
\int_X \alpha \wedge (w \vdash \beta).$$
\end{lemma}

\begin{proof}
Without loss of generality, we assume that the bi-vector field
 $w = w^{ij} e_i \wedge e_j$
for some vector fields $e_i$ over $X$, 
where $w^{ij}$ are smooth functions on $X$.
(Indeed, we use a partition of unity to decompose $w$
into a sum of bi-vector fields which can written this way.)
Notice that $(e_j \vdash \alpha) \wedge \beta$
and $\alpha \wedge (e_i \vdash \beta)$
have degree $|\alpha| + |\beta| -1 = \dim (X) + 1$,
hence they must  vanish.
Then we have
\begin{eqnarray*}
& & \int_X (w \vdash \alpha) \wedge \beta = 
\int_X w^{ij}(e_i \vdash e_j \vdash \alpha) \wedge \beta \\
& = & \int_X w^{ij} e_i \vdash [(e_j \vdash \alpha) \wedge \beta]
	- \int_X w^{ij} (-1)^{|\alpha|-1}
	(e_j \vdash \alpha) \wedge (e_i \vdash\beta) \\
& = & (-1)^{|\alpha|} 
\int_X w^{ij} e_j \vdash [\alpha \wedge (e_i \vdash \beta)]
- w^{ij} (-1)^{|\alpha|} \alpha 
\wedge (e_j \vdash e_i \vdash \beta) \\
& = & - \int_X w^{ij}\alpha 
\wedge (e_j \vdash e_i \vdash \beta)
= \int_X \alpha \wedge (w \vdash \beta).
\end{eqnarray*}
\end{proof}

\begin{proposition}
For any bi-vector field on a closed oriented $X$, we have
$$\int_X (\Delta \alpha) \wedge \beta = (-1)^{|\alpha|} 
\int_X \alpha \wedge \Delta \beta.$$
\end{proposition}

\begin{proof}
Using $\Delta \alpha = w \vdash d \alpha - d (w \vdash \alpha)$,
Lemma \ref{lem:contraction} and Stokes theorem, we have
\begin{eqnarray*}
&   & \int_X (\Delta \alpha) \wedge \beta = 
\int_X (w \vdash d \alpha) \wedge \beta 
- \int_X d(w \vdash \alpha) \wedge \beta \\
&  = & \int_X d \alpha \wedge (w \vdash \beta)
+(-1)^{|\alpha|} \int_X (w \vdash \alpha) \wedge d \beta \\
&  = & (-1)^{|\alpha|+1} \int_X \alpha \wedge d(w \vdash \beta)
+(-1)^{|\alpha|} \int_X\alpha \wedge (w \vdash d \beta) \\
& = & (-1)^{|\alpha|} \int_X \alpha \wedge \Delta \beta.
\end{eqnarray*}
\end{proof}

By Poincar\'{e} duality, $\int$ induces a nondegenerate pairing
on $H = H({\cal A} ,d )$, which is the de Rham cohomology.
Since $X$ is compact, $H$ is finite dimensional.
Thus only Condition $3$ in Theorem 1.1 remains to satisfy.
Thus, we have 

\begin{theorem} \label{thm:Poisson}
Assume that $(X, w)$ is a closed Poisson manifold 
such that the inclusions $i:(\Ker \Delta, d) \hookrightarrow
(\Omega(X), d)$ and 
$j:(\Ker d, \Delta) \hookrightarrow (\Omega(X), \Delta)$ 
induces isomorphisms  $H(i)$ and $H(j)$ on cohomologies.
If $K$ is the algebra of formal power series
on the de Rham cohomology ring of $X$,
then there is a structure of formal Frobenius manifold 
on the formal spectrum of $K$.
\end{theorem}

\begin{remark}
For a symplectic manifold $(X^{2n}, \omega)$,
Brylinski \cite{Bry} defined the symplectic star operator
$*_{\omega}: \Omega^k(X) \rightarrow \Omega^{2n-k}(X)$.
He showed that $\Delta = (-1)^k*_{\omega}d*_{\omega}$ on $\Omega^k(X)$
(Brylinski \cite{Bry}, Theorem 2.2.1,
where $\delta$ is used instead of $\Delta$).
From this it is clear that $H(i)$ is an isomorphism 
if and only if $H(j)$ is an isomorphism.
Hence it suffices to consider $H(i)$.
\end{remark}

On a symplectic manifold $(X^{2n}, \omega)$, notice that 
the inclusion 
$\phi: \Ker d  \cap \Ker \Delta \rightarrow (\Omega(X), d)$
factor through the inclusions $\psi: \Ker d  \cap \Ker \Delta
\rightarrow (\Ker \Delta, d)$ and
$i: (\Ker \Delta, d) \rightarrow (\Omega(X), d)$.
It is easy to see that $H(\psi)$ is surjective with kernel
$d\Ker \Delta$. Therefore, we have an isomorphism
$$\tilde{H}(\psi): \Ker d \cap \Ker \Delta / d\Ker \Delta 
	\cong H(\Ker \Delta, d).$$
Similarly, the kernel of $H(\phi)$ is $\Img d \cap \Ker \Delta$,
so we have an injective homomorphism
$$\tilde{H}(\phi): \Ker d \cap \Ker \Delta / \Img d \cap \Ker \Delta
\rightarrow H(X).$$
Also since $d \Ker \Delta \subset \Img d \cap \Ker \Delta$,
we have a surjective homomorphism 
$q:  \Ker d \cap \Ker \Delta / d\Ker \Delta 
\rightarrow \Ker d \cap \Ker \Delta / \Img d \cap \Ker \Delta$.
To summarize, we get a commutative diagram
$$\CD
\Ker d \cap \Ker \Delta / d\Ker \Delta
@>{\tilde{H}(\psi)}>> H(\Ker \Delta, d) \\
@VqVV @VV{H(i)}V \\
\Ker d \cap \Ker \Delta / \Img d \cap \Ker \Delta
@>>{\tilde{H}(\phi)}> H(X).
\endCD
$$
It is then clear that $H(i)$ is injective
if and only if $d \Ker \Delta = \Img d \cap \Ker \Delta$ (see also 
Manin \cite{Man}, (5.14)).
On the other hand, 
$H(i) = \tilde{H}(\phi) q \tilde{H}(\psi)^{-1}$,
from which we see that 
$H(i)$ is surjective if and only if  $\tilde{H}(\phi)$ is surjective.
The problem of surjectivity of $\tilde{H}(\phi)$  is 
equivalent to the following
question asked by Brylinski \cite{Bry}:
whether every class in $H(X)$ can  be represented 
by an element in $\Ker d \cap \Ker \Delta$.
He answered this question affirmatively for K\"{a}her manifolds.
For general symplectic manifolds, Mathieu \cite{Mat} and Yan \cite{Yan}
proved the following result by different methods:

\begin{proposition}
For any symplectic manifold $(X^{2n}, \omega)$, not necessarily compact,
the following two statements are equivalent:

(a). $\tilde{H}(\phi)$ is surjective.

(b). For each $0 \leq k \leq n$, $L^k: H^{n-k}(X) \rightarrow H^{n+k}(X)$
is surjective,
where $L$ is induced by wedge product with $\omega$.
\end{proposition}

For closed symplectic manifolds, Mathieu \cite{Mat} observed that 
(b) is equivalent to each $L^k$ being an isomorphism on $H^{n-k}(X)$
because $H^{n-k}(X)$ and $H^{n+k}(X)$ have the same dimension.
If (b) holds, one says that $(X, \omega)$ satisfies the hard Lefschetz
theorem. 
So $H(i)$ is surjective if and only if
the symplectic manifold $(X, \omega)$ satisfies the hard Lefschetz
theorem. 
One can find examples which do not satisfy the hard
Lefschetz theorem in Mathieu's paper.
So not every closed symplectic manifold satisfies the conditions in
Theorem \ref{thm:Poisson}.
Nevertheless, we will show that closed K\"{a}hler manifolds
do satisfy these conditions.

\subsection{de Rham Frobenius manifolds for K\"ahler manifolds}
In Cao-Zhou \cite{Cao-Zho2}, we constructed Frobenius manifold 
structure on the Dolbeault cohomology for K\"ahler manifolds. Here 
we carry out the construction on the de Rham cohomology.
 
For a K\"{a}hler manifold $(X, \omega)$,
let $L(\alpha) = \omega \wedge \alpha$,
and $\Lambda$ be the adjoint of the  operator $L$ 
defined by the Hermitian metric.
Then we have 
$$\Delta  = \Lambda d - \Lambda d = [\Lambda, d].$$ 
From the Hodge identities 
(Griffith-Harris \cite{Gri-Har}, p. 111),
we get
$$\Delta  = -4\pi (d^c)^* = 
\sqrt{-1}(\bar{\partial}^* - \partial^* ).$$
Where 
$$d^c = \frac{\sqrt{-1}}{4\pi}(\bar{\partial} - \partial).$$ 
It is clear that $\Delta^* = - 4\pi d^c
= - \sqrt{-1}(\bar{\partial} - \partial)$. 
Let $\square_{\Delta} = \Delta \Delta^* + \Delta^* \Delta$.

\begin{lemma} $\square_{\Delta} = \square$.
\end{lemma}

\begin{proof}
An easy calculation gives
\begin{eqnarray*}
\Delta^* \Delta & = & \partial \partial^*
	+ \bar{\partial}\bar{\partial}^* 
	- \bar{\partial} \partial^*
	- \partial \bar{\partial}^*, \\
\Delta \Delta^* & = & \partial^* \partial
	+ \bar{\partial}^* \bar{\partial}
	- \bar{\partial}^* \partial
	- \partial^* \bar{\partial}.
\end{eqnarray*}
Then we have
\begin{eqnarray*}
\square_{\Delta} = \square_{\partial} + \square_{\bar{\partial}}
+ ( \bar{\partial} \partial^* + \partial^* \bar{\partial}) 
-( \partial \bar{\partial^*} + \bar{\partial^*} \partial).
\end{eqnarray*}
The lemma follows from  the following well-known formulas in
K\"{a}hler geometry (Griffiths-Harris \cite{Gri-Har}, p. 115):
\begin{eqnarray} 
& \square = 2 \square_{\partial} 
= 2 \square_{\bar{\partial}}, \label{Laplacians} \\
& \bar{\partial} \partial^* + \partial^* \bar{\partial}
=  \partial \bar{\partial^*} + \bar{\partial^*} \partial = 0.
\end{eqnarray}
Here $\square =dd^*+d^*d$ is the Hodge Laplacian.
\end{proof}

Standard Hodge theory argument gives 
a decomposition 
$$\Omega(X) = {\cal H} \oplus \Img \Delta \oplus \Img \Delta^*,$$
where ${\cal H}$ is the space of harmonic forms on $X$.
Furthermore, the inclusion $({\cal H}, 0) \subset 
(\Omega(X), d)$ induces an isomorphism 
$H(\Omega(X), \Delta) \cong {\cal H}$.
Notice that $\Delta$ commutes with $\square_{\Delta}$,
so it commutes with $\square$.

\begin{lemma}\label{lemma:=0} On a K\"{a}hler manifold, 
$d \Delta^* + \Delta^* d= 0$
 and  $\Delta d^* + d^* \Delta= 0$.
\end{lemma}

\begin{proof} 
The first identity follows from the following identity
$$dd^c = - d^cd = \frac{\sqrt{-1}}{2\pi}\partial \bar{\partial}.$$
Taking formal adjoint gives the second identity.
\end{proof}

\begin{proposition} \label{prop:quasi}
Let $(X, \omega)$ be a closed K\"{a}hler manifold, 
then the  inclusions $i: (\Ker \Delta, d) \subset (\Omega(X), d)$ 
and $j: (\Ker d, \Delta) \subset (\Omega(X), \Delta)$ induce
isomorphisms on cohomology.
\end{proposition}

\begin{proof}
From Remark 2.1, it suffices to show that $H(i)$ is 
an isomorphism. 
By Hodge theory, every de Rham cohomology class on $X$ is represented 
by a harmonic form $\alpha$: $\square \alpha = 0$.
Now we have $\square =\square_{\Delta}$, 
so $\square_{\Delta}\alpha =0$, and hence  $\Delta \alpha = 0$. 
This shows that $H(i)$ is surjective. 

To show $H(i)$ is injective,
let $\alpha \in \Ker \Delta$, such that $\alpha = d\beta$ for some
$\beta \in \Omega(X)$. We need to show that $\alpha \in d \Ker \Delta$. 
By Hodge decomposition, $d\beta = \square\gamma$ 
for some $\gamma \in \Omega(X)$.
Hence we have, 
\begin{eqnarray*}
\square d \gamma = d \square \gamma = d(d\beta) = 0,
\end{eqnarray*}
which implies that $d\gamma = 0$. 
Similarly,
\begin{eqnarray*}
\square_{\Delta} \Delta\gamma = \Delta  \square_{\Delta} \gamma 
 = \Delta \alpha = 0
\end{eqnarray*}
implies that $\Delta\gamma = 0$. 
Now $\alpha = \square \gamma = 
(dd^* + d^*d)\gamma = d (d^*\gamma)$.
By Lemma \ref{lemma:=0},
$\Delta(d^*\gamma) = - d^* \Delta \gamma = 0$.
So $\alpha \in d \Ker \Delta$, 
 thus $H(i)$ is injective. 
\end{proof}

In conclusion, we have 

\begin{theorem} \label{thm:Frobenius}
For any closed K\"{a}hler manifold $X$,
if $K$ is the algebra of formal superpower series
on the de Rham cohomology ring of $X$,
then there is a structure of formal Frobenius manifold 
on the formal spectrum of $K$
obtained from Theorem \ref{thm:constrution} for the dGBV algebra
$(\Omega(X), \wedge, d, \Delta, [\cdot \bullet \cdot])$.
\end{theorem}

\subsection{Convergence of $\Gamma$} \label{sec:convergence}
For a closed K\"{a}hler  manifold $X$,
by Hodge theory, we can take $e_j$'s to be harmonic,
hence we automatically have $e_j \in \Ker d \cap \Ker \Delta$.
Furthermore,
\begin{eqnarray*}
\Img \Delta 
& = & \Delta ({\cal H} \oplus \Img d \oplus \Img d^*)\\
& = & \Img \Delta d \oplus \Img \Delta d^* \\
& = & \Img d \Delta \oplus \Img d^* \Delta.
\end{eqnarray*}
Now since $\Gamma_n\in \Img \Delta$ for $n \geq 2$,
we can further take $\Gamma_n$ to be in $\Img d^* \Delta$,
because we are inductively solving
$d\Gamma_n = \psi_n(\Gamma_1, \cdots, \Gamma_{n-1})$.

Let $G: \Omega(X) \rightarrow \Omega(X)$ be the Green's operator
of $\square$.

\begin{lemma} \label{lm:series}
Let $\Gamma = \sum_n \Gamma_n$ be a normalized universal 
solution to the Maurer-Cartan equation 
\begin{eqnarray} \label{eqn:dMC}
d \Gamma + \frac{1}{2}[\Gamma \bullet \Gamma] = 0
\end{eqnarray}
with 
$\Gamma_1 = \sum_j x^j e_j$, $e_j$ harmonic, 
$\Gamma_n \in \Img d^* \Delta$, for $n > 1$.
Then $\Gamma$ satisfies
$$\Gamma = \Gamma_1 - \frac{1}{2}d^*G[\Gamma \bullet \Gamma] = 0.$$
\end{lemma}

\begin{proof}
This is equivalent to solving the Maurer-Cartan equation 
inductively by imposing the above conditions.
\end{proof}

We call a solution as in Lemma \ref{lm:series} {\em analytically 
normalized}.
Now by modifying a standard argument in Kodaira-Spencer-Kuranishi
deformation theory (see e.g. Morrow-Kodaira \cite{Mor-Kod}.
Chapter 4, Proposition 2.4),
an analytically normalized solution $\Gamma$ is
convergent for small even $x^j$'s and all odd $x^j$'s.
(This method was also used by Tian \cite{Tia1} and Todorov \cite{Tod}.) 
Therefore, we actually obtain a Frobenius
supermanifold modelled on $H^*_{dR}(X)$,
with its Bosonic part a neighborhood of the zero vector in
$H^{\text{even}}_{dR}(X)$.

\section{Comparison with Frobenius supermanifold modelled on
Dolbeault cohomology} In this section, we identify the Frobenius 
manifold constructed in Theorem \ref{thm:Frobenius} with the one we 
constructed on the Dolbeault cohomology in \cite{Cao-Zho2}.

\subsection{Frobenius supermanifold modelled on
Dolbeault cohomology}
We review the construction of Cao-Zhou \cite{Cao-Zho2} 
in this section.
Let $(X, g, J)$ be a closed 
K\"{a}hler manifold with K\"{a}hler form $\omega$.
To furnish the comparison with Frobenius supermanifold 
structure on the de Rham cohomology in Theorem \ref{thm:Frobenius}, 
we make a slight modification.
Consider the quadruple 
$(\Omega^{*, *}(X), \wedge, \delta=\bar{\partial}, 
\Delta=-\sqrt{-1}\partial^*)$.
It is well-known that $\bar{\partial}^2 = 0$,
$(\partial^*)^2 = 0$, and 
$\bar{\partial}\partial^* + \partial^*\bar{\partial} = 0$.
Also,  $\bar{\partial}$ is a derivation.
Set 
$$[a \bullet b]_{\partial^*} = 
-\sqrt{-1}(-1)^{|a|}(\partial^*(a \wedge b ) 
- \partial^* a \wedge b - (-1)^{|a|}a \wedge \partial^* b).$$
Then $(\Omega^{*, *}(X), \wedge, \delta=\bar{\partial}, 
\Delta=-\sqrt{-1}\partial^*, 
[\cdot\bullet\cdot]_{\partial^*})$ is a dGBV algebra.
Furthermore, let $\int_X: \Omega^{*, *}(X) \rightarrow {\Bbb C}$
be the ordinary integration of differential forms.
Then $\int_X$ is a nice integral for the above dGBV algebra.
Hodge theoretical argument similar to the ones in
last section shows that the two natural inclusions 
$i: (\Ker \partial^*, \bar{\partial}) \rightarrow 
	(\Omega^{*, *}(X), \bar{\partial})$
and 
$j: (\Ker \bar{\partial},  \partial^*) \rightarrow 
	(\Omega^{*, *}(X), \partial^*)$
induce isomorphisms on cohomology.
Therefore, Theorem \ref{thm:constrution} applies.
To see that we get a Frobenius supermanifold this way,

we use, similar to what we did in \S \ref{sec:convergence}, the Hodge decomposition of 
$\square_{\bar{\partial}}$ to get
\begin{eqnarray*}
\Img \partial^* 
& = & \partial^*({\cal H} \oplus \Img \bar{\partial} \oplus
	\Img \bar{\partial}^*) \\
& = & \Img \partial^*\bar{\partial} \oplus 
	\Img  \partial^*\bar{\partial}^* \\
& = & \Img \bar{\partial}\partial^* \oplus 
	\Img  \bar{\partial}^*\partial^*.
\end{eqnarray*}
We have used the fact 
$\partial^*\bar{\partial} + \bar{\partial}\partial^* = 0$
and $\partial^*\bar{\partial}^* + \bar{\partial}^*\partial^* = 0$
above.
Similarly, we can find a  normalized universal  
solution $\Gamma = \sum_n \Gamma_n$ to 
\begin{eqnarray} \label{eqn:dbarMC}
\bar{\partial} \Gamma 
	+ \frac{1}{2}[\Gamma\bullet\Gamma]_{\partial^*}
	= 0,
\end{eqnarray}
such that $\Gamma_1 = \sum_j x^j e_j$, 
$e_j$ $\bar{\partial}$-harmonic, and
$\Gamma_n \in \Img \bar{\partial}^* \partial^*$, for $n > 1$.
We also call such a solution analytically normalized.
The convergence of such a $\Gamma$ can be established 
as in \S \ref{sec:convergence}.
This gives rise to a 
Frobenius supermanifold structure modelled on the Dolbeault
cohomology.

\begin{remark}
Similarly, set 
$$[a \bullet b]_{\bar{\partial}^*} = 
\sqrt{-1}(-1)^{|a|}(\bar{\partial}^*(a \wedge b ) 
- \bar{\partial}^* a \wedge b 
- (-1)^{|a|}a \wedge \bar{\partial}^* b).$$
Then $(\Omega^{*, *}(X), \wedge, \delta={\partial}, 
\Delta=\sqrt{-1}\bar{\partial}^*, 
[\cdot\bullet\cdot]_{\bar{\partial}^*})$ is a dGBV algebra.
\end{remark}

\subsection{ The identification}
By Hodge theory, there is a natural isomorphism between 
$H(\Omega^{*, *}(X), \bar{\partial})$ 
and ${\cal H}_{\bar{\partial}}$, the 
space of $\bar{\partial}$-harmonic forms. 
Since $\square_{\bar{\partial}} = \frac{1}{2}\square$,
we have
\begin{eqnarray*}
\square_{\bar{\partial}}\bar{\alpha} 
= \frac{1}{2}\square \bar{\alpha}
= \frac{1}{2}\overline{\square \alpha}
= \overline{\square_{\bar{\partial}}\alpha}.
\end{eqnarray*}
Hence, complex conjugation gives ${\cal H}_{\bar{\partial}}$
a real structure.
It follows then we can take $e_j$'s to be real. 
Therefore, with respect to the induced real structure on $K$,
$x^j$'s are also real.

\begin{lemma}
There is a unique analytically normalized solution $\Gamma$ to
(\ref{eqn:dbarMC}) 
which also satisfies
\begin{eqnarray} \label{eqn:rondMC}
\partial \Gamma 
+ \frac{1}{2}[\Gamma\bullet\Gamma]_{\bar{\partial}^*} = 0.
\end{eqnarray}
Furthermore, $\Gamma$ is real, and is an
analytically normalized solution to (\ref{eqn:dMC}).
\end{lemma}

\begin{proof}
This can be proved inductively.
By Hodge theory  for $\bar{\partial}$ and $\partial$, 
we have the following decompositions
\begin{eqnarray*}
\Omega^{*,*}(X) 
 =  {\cal H} \oplus \Img \bar{\partial} 
	\oplus \Img \bar{\partial}^* 
 =  {\cal H} \oplus \Img \partial \oplus \Img \partial^*.
\end{eqnarray*}
Use these decomposition twice, we get
$$\Omega^{*,*}(X) 
 = {\cal H} \oplus \Img \bar{\partial} \partial
	\oplus \Img \bar{\partial}^* \partial
 \oplus \Img \bar{\partial} \partial^*
	\oplus \Img \bar{\partial}^* \partial^*.$$
Denote the Green's operators for $\square$, 
$\square_{\bar{\partial}}$ and $\square_{\partial}$ by
$G$, $G_{\bar{\partial}}$ and $G_{\partial}$ respectively.
Then we have $G_{\bar{\partial}} = G_{\partial} = 2G$.
Since $\Gamma_1 = \sum_j x^j e_j$, $e_j$'s harmonic,
we have
\begin{eqnarray*}
[\Gamma_1\bullet\Gamma_1]_{\partial^*} 
& = & -\sqrt{-1}(\partial^*(\Gamma_1 \wedge\Gamma_1) 
- (\partial^*\Gamma_1) \wedge \Gamma_1 
- \Gamma_1 \wedge \partial^*(\Gamma_1) ) \\
& = & -\sqrt{-1} \partial^*(\Gamma_1 \wedge\Gamma_1).
\end{eqnarray*}
Similarly we have
$$[\Gamma_1\bullet\Gamma_1]_{\bar{\partial}^*}
=\sqrt{-1} \bar{\partial}^*(\Gamma_1 \wedge\Gamma_1).$$
So we need to simultaneously solve
\begin{align*}
\bar{\partial} \Gamma_2 & = \frac{1}{2}\sqrt{-1}
\partial^*(\Gamma_1 \wedge\Gamma_1), \tag{$\bar{\partial}_2$} \\
\partial \Gamma_2 & = -\frac{1}{2}\sqrt{-1}
\bar{\partial}^*(\Gamma_1 \wedge\Gamma_1), \tag{$\partial_2$}
\end{align*}
for $\Gamma_2 \in \Img \bar{\partial}^*\partial^*$.
Taking $\bar{\partial}^*$ on both sides of 
($\bar{\partial}_2$), we get
$$\square_{\bar{\partial}}\Gamma_2 
= \bar{\partial}^*\bar{\partial} \Gamma_2
= \frac{1}{2}\sqrt{-1}
\bar{\partial}^*\partial^*(\Gamma_1 \wedge\Gamma_1).$$
Therefore, 
$$\Gamma_2 = \frac{1}{2}\sqrt{-1} G_{\bar{\partial}}
\bar{\partial}^*\partial^*(\Gamma_1 \wedge\Gamma_1)
= -\frac{1}{2}\sqrt{-1} G_{\partial}
 \partial^*\bar{\partial}^*(\Gamma_1 \wedge\Gamma_1)$$
is the unique candidate for the  solution. 
Now we have
\begin{eqnarray*}
&   &   \bar{\partial} \Gamma_2 = 
 \frac{1}{2}\sqrt{-1} \bar{\partial} G_{\bar{\partial}}
 \bar{\partial}^*\partial^*(\Gamma_1 \wedge\Gamma_1) \\
& = & \frac{1}{2}\sqrt{-1}  G_{\bar{\partial}} (\bar{\partial}
 \bar{\partial}^*)\partial^*(\Gamma_1 \wedge\Gamma_1) 
 = \frac{1}{2}\sqrt{-1} G_{\bar{\partial}}
(\square_{\bar{\partial}} - \bar{\partial}^*\bar{\partial}) 
 \partial^*(\Gamma_1 \wedge\Gamma_1) \\
& = &  \frac{1}{2}\sqrt{-1}  G_{\partial}
\square_{\partial} \partial^*(\Gamma_1 \wedge\Gamma_1) 
 =   \frac{1}{2}\sqrt{-1} 
	\partial^* (\Gamma_1 \wedge\Gamma_1) \\
& = & - \frac{1}{2} [\Gamma_1\bullet\Gamma_1]_{{\partial}^*}.
\end{eqnarray*}
So $\Gamma_2$ satisfies $(\bar{\partial}_2)$.
Similarly, one can check that $\Gamma_2$ satisfies 
$(\partial_2)$.
The general induction procedure is similar.
Now $\overline{\Gamma}$ is also an analytically normalized
solution to both (\ref{eqn:dbarMC}) and (\ref{eqn:rondMC}),
we have $\overline{\Gamma} = \Gamma$.
Furthermore, 
adding (\ref{eqn:dbarMC}) to (\ref{eqn:rondMC}) shows
that $\Gamma$ also satisfies (\ref{eqn:dMC}).
\end{proof}

As a corollary, we get

\begin{theorem}
For a closed K\"{a}hler manifold $X$,
there is a real structure on the Frobenius supermanifold 
modelled on Dolbeault cohomology obtained from Theorem \ref{thm:constrution}.
Its real part can be identified with the Frobenius supermanifold 
modelled on de Rham cohomology obtained from Theorem \ref{thm:constrution}.
\end{theorem}

\end{document}